\documentclass[a4paper,11pt,twoside,reqno]{amsart}
\usepackage{latexsym,amsmath,amsthm,verbatim,ifthen,amssymb,graphicx,color,mathrsfs}
\usepackage[all]{xy}
\usepackage{color}
\usepackage{hyperref}
\hypersetup{
    colorlinks,
    citecolor=blue,
    filecolor=black,
    linkcolor=blue,
    urlcolor=black
}
\usepackage[a4paper]{geometry}
 \title[Spatial realization and the Bar construction]{Spatial realization of a Lie algebra and the Bar construction}

\author{Yves F\'elix}
\address{Institut de Math\'ematiques et Physique\\
         Universit\'e Catholique de Louvain-la-Neuve\\
         Louvain-la-Neuve\\
         Belgique}
\email{Yves.felix@uclouvain.be}

\author{Daniel Tanr\'e}
\address{D\'epartement de Math{\'e}matiques\\
         UMR 8524 \\
         Universit\'e de Lille\\
         59655 Villeneuve d'Ascq Cedex\\
         France}
\email{Daniel.Tanre@univ-lille.fr}

\thanks{The authors are partially supported by the 
MINECO-FEDER grant  MTM2016-78647-P.
The second author is partially supported by the  ANR-11-LABX-0007-01  ``CEMPI''}

\subjclass[2010]{Primary: 55P62, 17B55; Secondary: 55U10}

\keywords{Rational homotopy theory. Realization of Lie algebras. Lie models of simplicial sets. Model category.}

\SelectTips{cm}{}
\newtheorem{theoremv}{Main Theorem}

 \theoremstyle{definition}
 
 \theoremstyle{remark}

 \numberwithin{equation}{section}

\def\L{\mathbb{L}}

\def\hL{{\widehat{\mathbb L}}}

\def\ad{{\rm ad}}

\def\lasu{{\mathfrak{L}}}
\def\cdgl{{\mathrm{cdgl}}}
\def\sset{{\mathbf{Sset}}}
\def\CDGL{{\mathbf{cdgl}}}
\def\Hom{{\mathrm{Hom}}}
\newcommand{\secref}[1]{Section~\ref{#1}}

\begin{document}

\date{\today}

\maketitle

\begin{abstract}
We prove that the spatial realization of a rational complete Lie algebra $L$, concentrated in degree~0, 
is isomorphic to the 
simplicial bar construction on the group, obtained from the Baker-Campbell-Hausdorff product on 
 $L$.
\end{abstract}

\section*{Introduction}
In \cite{four}, we construct a cosimplicial differential graded complete Lie algebra  
(henceforth $\cdgl$) $(\lasu_{\bullet},d)$,
  in which $(\lasu_{1},d)$ is the Lawrence-Sullivan model of the interval introduced in \cite{LS}.
  As in the work of Sullivan (\cite{Dennis}) for differential commutative graded algebras, the existence of this
  cosimplicial object gives an adjoint pair of functors between the category $\CDGL$ of cdgl's and the
  category $\sset$ of simplicial sets, see \cite{four} or \cite{four3}.
  In this work, we focus on one of them, the spatial realization functor,
  $$\langle - \rangle \colon \CDGL \to \sset,$$
  defined by $\langle L\rangle =\Hom_{\CDGL}(\lasu_{\bullet},L)$ for $L\in\CDGL$.
  (Let us also notice that  $\langle L\rangle$ is isomorphic to the nerve of $L$,
 a deformation retract of the Getzler-Hinich realization, see \cite{four6}, \cite{Nicoud}.)
  
  \smallskip
 More precisely, we are interested in the realization $\langle L\rangle$ of a complete Lie algebra, $L$,
 concentrated in degree 0 and (thus) with the differential 0. In this case, a group structure can be defined on
 the set $L$ from the Baker-Campbell-Hausdorff formula. We denote by $\exp\;L$ this group.
 The realization $\langle L \rangle$ is an Eilenberg-MacLane space of type $K(\pi,1)$, see \cite{MR2440265}.
 The purpose of this work is the determination of $\langle L \rangle$ up to isomorphism.
 
 \begin{theoremv}
 Let $L$ be a complete differential graded Lie algebra, 
 concentrated in degree 0. 
 Then,  its spatial realization $\langle L\rangle $ is isomorphic 
 to the simplicial bar construction on $\exp\;L$.
 \end{theoremv}
 
 Let us emphasize that for an Eilenberg-Maclane space $K(G,1)$, of associated Lie algebra $L$, then the 
 realizations of the  $A_{PL}(K(G,1))$ of Sullivan (\cite{Dennis}) and ${\mathrm{MC}}_{*}(L)$ of Getzler 
 (\cite{Ge}) are 
 weakly equivalent to the simplicial bar resolution. In our setting, the simplicial set $\Hom({\mathcal L}_{\bullet},L)$ is
 isomorphic to the simplicial bar resolution. With \cite{four3}, we know that   spaces more general than a $K(G,1)$ admits 
 a cdgl model $L$. For them,  the simplicial set
 $\Hom({\mathcal L}_{\bullet},L)$  appears as a natural extension of the bar construction. We will come
 back on this point in a forthcoming work.
 
 \medskip
 In \secref{sec:recall}, we recall basic background on $\CDGL$ and our construction $\lasu_{\bullet}$. 
 \secref{sec:proof} consists of the proof of the Main theorem.
  
  \section{Some reminders}\label{sec:recall}
  We first recall the construction of the cosimplicial $\cdgl$ $\lasu_{\bullet}$.
 Let  $V$ be a finite dimensional graded vector space. 
 The completion of the free graded Lie algebra on $V$, $\L(V)$,
 is the inverse limit,
 $$\hL(V) = \varprojlim_n \L(V)/ \mathbb L^{\geq n} (V),$$
 where $\L^{\geq n} (V)$ is the ideal generated by the Lie brackets of length $\geq n$. 
 We call $\hL(V)$ the free complete graded Lie algebra on $V$.
 
 \medskip
As a graded Lie algebra, $\lasu_n$ is the  free complete graded Lie algebra 
on the rational vector space generated by the elements 
$a_{i_0\dots i_k}$  of degree $\vert a_{i_0\dots i_k}\vert = k-1$, with $0\leq i_0<\dots <i_k\leq n$.
We denote by $\ad_{i_0\dots i_k}$ the Lie derivation  $[a_{i_0\dots i_k},-]$. 
The $\cdgl$ $\lasu_{n}$ satisfies the following properties.
\\
-- $\lasu_0$ is the free $\cdgl$ on a Maurer-Cartan element $a_0$, that is:
$$\lasu_0 = (\L(a_0),d)\,, \quad da_0 = -\frac{1}{2} [a_0, a_0].$$
-- $\lasu_1 = (\hL(a_0, a_1, a_{01}),d)$ is the Lawrence-Sullivan interval (see \cite{LS}),
where $a_{0}$ and $a_{1}$ are Maurer-Cartan elements and
$$da_{01}=[a_{01},a_{1}]+\frac{\ad_{01}}{e^{\ad_{01}}-1}(a_{1}-a_{0}).$$
--  A model $\lasu_2$ for the triangle has been described in \cite{four} (see also \cite{four7}):
$$\lasu_2 = (\hL(a_0,a_1,a_2, a_{01}, a_{02}, a_{12}, a_{012}), d)
\text{ with }d(a_{012} ) = a_{01}* a_{12}*a_{02}^{-1}-[a_{0},a_{012}].$$
Here $*$ denotes the Baker-Campbell-Hausdorff product defined for any pair of elements $a$, $b$ of degree 0 by
$a\ast b=\log(\exp^a \,\exp^b)$.\\
-- Moreover these structures appear in each $\lasu_n$: 
each vertex $a_r$ is a Maurer-Cartan element, 
each triple $(a_r, a_s, a_{rs})$ is a Lawrence-Sullivan interval and
 each family 
 \linebreak
 $(a_r,a_s,a_t, a_{rs}, a_{rt}, a_{st}, a_{rst})$ is a triangle as above.

\medskip
The family $(\lasu_n)_{n\geq 0}$ forms a cosimplicial $\cdgl$ which allows the definition of the spatial realization 
of $L\in\CDGL$ by,
$$\langle L\rangle := \Hom_{\CDGL} (\lasu_{\bullet}, L).$$
The cofaces $\delta^i$ and the codegeneracies $\sigma^i$ of the cosimplicial $\cdgl$ $\lasu_{\bullet}$ are defined by
\begin{equation}\label{equa:delta}
\delta^i a_{i_{0}\dots i_{p}}= a_{j_{0}\dots j_{p}}
\quad
\text{with}
\quad
j_{k}=\left\{
\begin{array}{ll}
i_{k}&\text{if } i_{k}<i,\\
i_{k}+1&\text{if } i_{k}\geq i,
\end{array}\right.
\end{equation}
\begin{eqnarray}
\sigma^i a_{i_{0}\dots i_{p}}&=&
\left\{
\begin{array}{ll}
0
\quad
\text{if}
\quad
\left\{i,i+1\right\}\subset \left\{i_{0},\dots,i_{p}\right\},
\label{equa:delta1}\\[.2cm]
a_{j_{0}\dots j_{p}}
\quad
\text{otherwise, with}
\quad
j_{k}=\left\{
\begin{array}{ll}
i_{k}&\text{if } i_{k}\leq i,\\
i_{k}-1&\text{if } i_{k}> i.
\end{array}\right.
\end{array}\right.
\end{eqnarray}

\section{Proof of the main theorem}\label{sec:proof}
The simplicial bar construction  on a group $G$, 
is the simplicial set $B_{\bullet}G$ with set of $n$-simplices $B_{n}G=G^n$. Its elements are denoted
$[g_1\vert \dots \vert g_n]$, with $g_i\in G$. 
 The  faces $d_i$ and degeneracies $s_i$ of $B_{\bullet}G$ are defined as follows:
 \begin{eqnarray}
  d_0[g_1\vert \dots \vert g_n]&=& [g_2\vert \dots \vert g_n],\\
 d_i[g_1\vert \dots \vert g_n] &=& [g_1\vert \dots \vert g_ig_{i+1}\vert \dots \vert g_n],
  \quad \text{for}\quad 0<i<n,\nonumber\\
  d_n[g_1\vert \dots \vert g_n] &=& [g_1\vert \dots \vert g_{n-1}].\nonumber
  \end{eqnarray}
  The degeneracy $s_i$ inserts the identity $e$ of $G$ in position $i$.
  
\medskip
Let $L\in\CDGL$ be generated in degree 0 and $f \colon \lasu_n \to L$ a morphism in $\CDGL$.
 For degree reasons, we have $f(a_{i_0\dots i_k}) = 0$ if $k\neq 1$. Moreover, since $f$ commutes the differential,
 from the definition of the differential in $\lasu_{2}$, we get
$$0 = df(a_{0rs}) = f(a_{0r})\ast  f(a_{rs})* f(a_{0s})^{-1}.$$
Therefore, for any $r,s>0$, we have 
$$f(a_{rs}) = f(a_{0r})^{-1} * f(a_{0s}).$$
The map $f$ being entirely defined by its values on the $a_{0i}$, we have a bijection
$$\Phi \colon \Hom_{\CDGL}(\lasu_n, L) \to L^n ,\quad \text{defined by}\quad
f\mapsto (f(a_{01}), f(a_{02}), \dots , f(a_{0n})).$$
We  now determine the image of the faces and degeneracies on $\Hom_{\CDGL}(\lasu_{\bullet},L)$,
induced from (\ref{equa:delta}) and (\ref{equa:delta1}).
For the face operators, as only the elements $(a_{0r})$ play a role,  it suffices to consider,
$$ \delta^i (a_{0r}) = \left\{ \begin{array}{ll} a_{0,r} &  \text{if}\;r<i\\ a_{0, r+1}&
 \text{if}\;r\geq i\end{array}\right.
\quad \text{for}\quad i>0,\quad \text{and}\quad \delta^0(a_{0,r}) = a_{1, r+1}.
$$
Let $f\colon\lasu_{n}\to L$ be specified by 
$(f(a_{01}),\dots,f(a_{0n}))$.
Then $d_{0}f=f\circ \delta^0\colon \lasu_{n-1}\to L$ is described by
\begin{eqnarray*}
d_{0}(f(a_{01}),\dots,f(a_{0(n-1)}))
&=&
(f\circ \delta^0(a_{01}),\dots,f\circ \delta^0(a_{0(n-1)}))\\
&=&
(f(a_{12}),\dots,f(a_{1n}))\\
&=&
(f(a_{01})^{-1}\ast f(a_{02}),\dots,f(a_{01})^{-1}\ast f(a_{0n})).
\end{eqnarray*}
Therefore, the face operator $d_{0}$ on $L^{\bullet}$, induced from $\Phi$, is
$$d_{0}(x_{1},\dots,x_{n})=(x_{1}^{-1}\ast x_{2},\dots,x_{1}^{-1}\ast x_{n}).$$
Similar arguments give, for $i>0$,
$$d_{i}(x_{1},\dots,x_{n})=(x_{1},\dots,\hat{x}_{i},\dots,x_{n}).$$
As for the degeneracies, starting from
$$\sigma^i( a_{0r}) = \left\{ \begin{array}{ll}
 a_{0r} & \text{if}\;r\leq i\\
 a_{0,r-1} &  \text{if}\;r>i\end{array}\right.
\quad \text{for}\quad i>0,$$
$$\sigma^0(a_{0r}) = a_{0, r-1}\quad \text{if} \quad  r>1,
\quad \text{and}\quad
\sigma^0(a_{0,1}) =0,
$$
we get
 $$s_0(x_1, \dots , x_n) = (0, x_1, \dots , x_n)$$
and for $i>0$,
$$s_i (x_1, \dots , x_n) = (x_1, \dots , x_i, x_i, \dots , x_n).$$
Now a straightforward and easy computation shows that the morphism
$$\Psi \colon \Hom_{\CDGL}(\lasu_{\bullet}, L) \to B_{\bullet}(\exp L)$$ defined by 
$$\Psi(f)=[f(a_{01})|f(a_{01})^{-1}f(a_{02})|f(a_{02})^{-1}f(a_{03})|\dots|f(a_{0(n-1)})^{-1}f(a_{0n})]$$
is an isomorphism of simplicial sets.
\providecommand{\bysame}{\leavevmode\hbox to3em{\hrulefill}\thinspace}
\providecommand{\MR}{\relax\ifhmode\unskip\space\fi MR }
\providecommand{\MRhref}[2]{%
  \href{http://www.ams.org/mathscinet-getitem?mr=#1}{#2}
}
\providecommand{\href}[2]{#2}


\begin{thebibliography}{1}

\bibitem{four6}
Urtzi {Buijs}, Yves {F\'elix}, Aniceto {Murillo}, and Daniel {Tanr\'e},
  \emph{{The infinity Quillen functor, Maurer-Cartan elements and DGL
  realizations}}, arXiv e-prints (2017), arXiv:1702.04397.

\bibitem{four3}
\bysame, \emph{Homotopy theory of complete {L}ie algebras and {L}ie models of
  simplicial sets}, J. Topol. \textbf{11} (2018), no.~3, 799--825. \MR{3989431}

\bibitem{four7}
\bysame, \emph{Symmetric {L}ie models of a triangle}, Fund. Math. \textbf{246}
  (2019), no.~3, 289--300. \MR{3959254}

\bibitem{four}
\bysame, \emph{Lie models of simplicial sets and representability of the
  {Q}uillen functor}, Israel J. Math. \textbf{238} (2020), no.~1, 313--358.
  \MR{4145802}

\bibitem{MR2440265}
Xue~Zhi Cheng and Ezra Getzler, \emph{Transferring homotopy commutative
  algebraic structures}, J. Pure Appl. Algebra \textbf{212} (2008), no.~11,
  2535--2542. \MR{2440265}

\bibitem{Ge}
Ezra Getzler, \emph{Lie theory for nilpotent {$L_\infty$}-algebras}, Ann. of
  Math. (2) \textbf{170} (2009), no.~1, 271--301. \MR{2521116}

\bibitem{LS}
Ruth Lawrence and Dennis Sullivan, \emph{A formula for topology/deformations
  and its significance}, Fund. Math. \textbf{225} (2014), no.~1, 229--242.
  \MR{3205571}

\bibitem{Nicoud}
Daniel Robert-Nicoud, \emph{Representing the deformation {$\infty$}--groupoid},
  Algebr. Geom. Topol. \textbf{19} (2019), no.~3, 1453--1476. \MR{3954288}

\bibitem{Dennis}
Dennis Sullivan, \emph{Infinitesimal computations in topology}, Inst. Hautes
  \'Etudes Sci. Publ. Math. (1977), no.~47, 269--331 (1978). \MR{0646078 (58
  \#31119)}

\end{thebibliography}
\end{document}